\documentclass[a4paper,12pt]{amsart}
\usepackage{amssymb,  amsfonts, amsthm}
\usepackage{ifthen}
\usepackage{graphicx}
\nonstopmode \numberwithin{equation}{section}
\usepackage[usenames]{color}
\newtheorem{thm}{Theorem}
\newtheorem{cor}{Corollary}
\newtheorem{lem}{Lemma}
\newtheorem{prop}{Proposition}
\newtheorem{rem}{Remark}[section]
\newtheorem{rems}[equation]{Remarks}
\theoremstyle{definition}
\newtheorem{defin}{Definition}
\newtheorem{examp}[equation]{Example}
\newtheorem{prob}[equation]{Problem}
\newtheorem{ques}[equation]{Question}
\newtheorem{op}{Open Problem}

\newtheorem{conj}[equation]{Conjecture}
\newtheorem{deter}[equation]{Determination}
\newtheorem{case}{Case}[section]
\newtheorem{subcase}[equation]{Subcase}
\newtheorem{claim}{Claim}[section]
\newtheorem{subclaim}{Subclaim}

\newcounter {own}
\def\theown {\thesection       .\arabic{own}}

\newenvironment{pf}[1][]{%
 \vskip 3mm
 \noindent
 \ifthenelse{\equal{#1}{}}%
  {{\bf Proof. }}%
  {{\bf #1.} }%
 }%
{\qed\bigskip}

\newcounter{alphabet}
\newcounter{tmp}
\newenvironment{Thm}[1][]{\refstepcounter{alphabet}%
\bigskip%
\noindent%
{\bf Theorem \Alph{alphabet}}%
\ifthenelse{\equal{#1}{}}{}{ (#1)}%
{\bf .} \itshape}{\vskip 8pt}
\makeatletter
\newcommand{\Ref}[1]{\@ifundefined{r@#1}{}{\setcounter{tmp}{\ref{#1}}\Alph{tmp}}}
\makeatother

\newcommand{\IR}{{\mathbb R}}

\newcommand{\diam}{{\operatorname{diam}}}

\def\be{\begin{equation}}
\def\ee{\end{equation}}

\newcommand{\bee}{\begin{enumerate}}
\newcommand{\eee}{\end{enumerate}}

\newcommand{\blem}{\begin{lem}}
\newcommand{\elem}{\end{lem}}
\newcommand{\bthm}{\begin{thm}}
\newcommand{\ethm}{\end{thm}}
\newcommand{\bcor}{\begin{cor}}
\newcommand{\ecor}{\end{cor}}
\newcommand{\beg}{\begin{examp}}
\newcommand{\eeg}{\end{examp}}
\newcommand{\begs}{\begin{examples}}
\newcommand{\eegs}{\end{examples}}
\newcommand{\bdefe}{\begin{defin}}
\newcommand{\edefe}{\end{defin}}
\newcommand{\bprob}{\begin{prob}}
\newcommand{\eprob}{\end{prob}}
\newcommand{\bques}{\begin{ques}}
\newcommand{\eques}{\end{ques}}
\newcommand{\bei}{\begin{itemize}}
\newcommand{\eei}{\end{itemize}}

\newcommand{\bde}{\begin{deter}}
\newcommand{\ede}{\end{deter}}
\newcommand{\bca}{\begin{case}}
\newcommand{\eca}{\end{case}}
\newcommand{\bsca}{\begin{subcase}}
\newcommand{\esca}{\end{subcase}}
\newcommand{\bcl}{\begin{claim}}
\newcommand{\ecl}{\end{claim}}
\newcommand{\bscl}{\begin{subclaim}}
\newcommand{\escl}{\end{subclaim}}
\newcommand{\bcon}{\begin{conj}}
\newcommand{\econ}{\end{conj}}
\newcommand{\bcons}{\begin{conjs}}
\newcommand{\econs}{\end{conjs}}
\newcommand{\bprop}{\begin{propo}}
\newcommand{\eprop}{\end{propo}}
\newcommand{\br}{\begin{rem}}
\newcommand{\er}{\end{rem}}
\newcommand{\brs}{\begin{rems}}
\newcommand{\ers}{\end{rems}}
\newcommand{\bo}{\begin{obser}}
\newcommand{\eo}{\end{obser}}
\newcommand{\bos}{\begin{obsers}}
\newcommand{\eos}{\end{obsers}}
\newcommand{\bpf}{\begin{pf}}
\newcommand{\epf}{\end{pf}}
\newcommand{\ba}{\begin{array}}
\newcommand{\ea}{\end{array}}
\newcommand{\beq}{\begin{eqnarray}}
\newcommand{\beqq}{\begin{eqnarray*}}
\newcommand{\eeq}{\end{eqnarray}}
\newcommand{\eeqq}{\end{eqnarray*}}

\newcommand{\ds}{\displaystyle}

\newcommand{\bop}{\begin{op}}
\newcommand{\eop}{\end{op}}

\newtheorem{pfofThm1.5}[equation]{}

\newcounter{minutes}
\divide\time by 60
\newcounter{hours}
\multiply\time by 60 \addtocounter{minutes}{-\time}

\begin{document}

\bibliographystyle{amsplain}

\title{On the subinvariance of uniform domains in Banach spaces}

\def\thefootnote{}
\footnotetext{ \texttt{\tiny File:~\jobname .tex,
          printed: \number\year-\number\month-\number\day,
          \thehours.\ifnum\theminutes<10{0}\fi\theminutes}
} \makeatletter\def\thefootnote{\@arabic\c@footnote}\makeatother

\author{M. Huang}
\address{M. Huang, Department of Mathematics,
Hunan Normal University, Changsha,  Hunan 410081, People's Republic
of China} \email{mzhuang79@yahoo.com.cn}

\author{M. Vuorinen}
\address{M. Vuorinen, Department of Mathematics and Statistics,
University of Turku, 20014 Turku, Finland } \email{vuorinen@utu.fi}

\author{X. Wang $^* $
}
\address{X. Wang, Department of Mathematics,
Hunan Normal University, Changsha,  Hunan 410081, People's Republic
of China} \email{xtwang@hunnu.edu.cn}

\date{}
\keywords{ Uniform domain, quasiconformal mapping,  FQC mapping, QH homeomorphism, (c,h)-solid curve.\\
${}^{\mathbf{*}}$ Corresponding author}

\begin{abstract}
Suppose that $E$ and $E'$ denote real Banach spaces with dimension at least $2$, that $D\subset E$ and $D'\subset E'$ are
domains, and that $f: D\to D'$ is a homeomorphism. In this paper,  we prove the following
subinvariance property for the class of uniform domains: Suppose that $f$ is a freely quasiconformal
mapping and that $D'$ is uniform. Then the image $f(D_1)$ of every uniform subdomain $D_1$ in $D$ under $f$ is still uniform.
This result answers an open problem of V\"ais\"al\"a
 in the affirmative.
\end{abstract}

\thanks{The research was partly supported by NSF of
China (No. 11071063 and 11101138).
The research of Matti Vuorinen was supported by the Academy of Finland,
Project 2600066611. }

\maketitle\pagestyle{myheadings} \markboth{}{On the subinvariance of
uniform domains in Banach spaces}

\section{Introduction and main results}\label{sec-1}

The theory of quasiconformal mappings in Banach spaces was developed
by V\"ais\"al\"a in a series of papers \cite{Vai6-0}$-$\cite{Vai8}
published in the period 1990$-$1992. In infinite dimensional cases
the classical methods such as the extremal length and conformal
invariants are no longer available. V\"ais\"al\"a built his theory
using notions from the theory of metric spaces. The basic notions
are curves and their lengths as well as special classes of domains
such as uniform and John domains. In fact, a given domain $D$ in a
Banach space $E$ it is essential to consider several metric space
structures, including hyperbolic type metrics of $D,$ at the same
time. The basic metrics are the norm metric of $E$ and the
quasihyperbolic and distance ratio metrics of $D$.

It is a natural question to study which properties of quasiconformal
mappings in the Euclidean spaces have their counterparts for freely
quasiconformal mappings of Banach spaces. The term "free" in this
context was coined by V\"ais\"al\"a and it refers to the
dimensionfree character of the results in the case of infinite
dimensional Banach spaces. Due to V\"ais\"al\"a's work, many results
of this type are already known. In the same spirit, we will
investigate a subinvariance problem.

It is well known that uniform domains are subinvariant under quasiconformal
mappings in $\IR^n$. By this, we mean that if $f:\; D\to D'$ is a $K$-quasiconformal
mapping between domains in $\IR^n$ and the image domain $D'$ is $c$-uniform,
then $D_1'=f(D_1)$ is $c'$-uniform,
 for every $c$-uniform subdomain
 $D_1 \subset D$,
 where $c'=c'(c, K, n)$. This follows from the corresponding result for QED domains
\cite[Remark]{FHM} and from \cite[Theorem 5.6]{Vai2}.

Our work is motivated by these ideas which we will extend to the
context of freely quasiconformal mappings in Banach spaces.   In his
study of the subinvariance of uniform domains in Banach spaces
V\"ais\"al\"a proved the following invariance property under QH
mappings (see Definition \ref{def1.6} in Section \ref{sec-2}); see
\cite[Theorem 2.44]{Vai4}.

\begin{Thm}\label{ThmB}
Suppose that $D\subset E$ and $D'\subset E'$, where $E$ and $E'$ are Banach spaces with dimension at least $2$,
that the domains $D'$ and $G\subset D$ are $c$-uniform, and that $f:\; D\to D'$ is $M$-QH.
Then $f(G)$ is $c'$-uniform with $c'=c'(c,M)$.
\end{Thm}

A natural problem is  whether the corresponding result holds for
freely quasiconformal mappings (see Definition \ref{def1.7} in
Section \ref{sec-2}) in Banach spaces or not. In fact, the problem
was raised by V\"ais\"al\"a in the following way; see
\cite[Subsection 2.43]{Vai4}.

\bop\label{Con1} Suppose that $D'\subset E'$ is $a$-uniform and that
$f:D\to D'$ is a $\psi$-FQC mapping (i.e., a freely
$\psi$-quasiconformal mapping), where $D\subset E$. If
$D_1$ is a $c$-uniform subdomain  in $D$, is it true that then $D_1'=f(D_1)$ is $c'$-uniform with $c'=c'(a, c , \psi)$? \eop

Our main result is the following theorem, which shows that the answer to Open
Problem \ref{Con1} is in the affirmative.

\begin{thm}\label{thm1.1}
Suppose that $D'\subset E'$ is $a$-uniform and that $f:D\to D'$ is a $\psi$-FQC mapping, where $D\subset E$. For each
subdomain $D_1$ in $D$, if $D_1$ is $c$-uniform, then $f(D_1)$ is still $c'$-uniform, where $c'=c'(a, c , \psi)$.
\end{thm}

We remark that the hypothesis ``$D'$ being uniform" in Theorem
\ref{thm1.1} is necessary. This can be seen by letting
$D=\mathbb{B}$ and $D'=\mathbb{B}\backslash [0,1)$ in $\IR^2$.

The proof of Theorem \ref{thm1.1} will be given
in Section \ref{sec-4}.
 In Section \ref{sec-2}, some necessary preliminaries will be introduced.

\section{Preliminaries}\label{sec-2}
\subsection{Notation}
Throughout the paper, we always assume that $E$ and $E'$ denote real
Banach spaces with dimension at least $2$. The norm of a vector $z$
in $E$ is written as $|z|$, and for every pair of points $z_1$,
$z_2$ in $E$, the distance between them is denoted by $|z_1-z_2|$,
the closed line segment with endpoints $z_1$ and $z_2$ by $[z_1,
z_2]$. Moreover, we use $\mathbb{B}(x, r)$ to denote the ball with center $x\in E$ and radius $r$ $(\geq 0)$,
 and its boundary and closure are denoted by
$\mathbb{S}(x,\; r)$ and $\overline{\mathbb{B}}(x,\; r)$,
respectively. In particular, we use $\mathbb{B}$ to denote the unit
ball $\mathbb{B}(0,\; 1)$.

For the convenience of notation, given domains $D \subset E\,,$ $D'
\subset E'$  and a mapping $f:D \to D'$ and points $x$, $y$, $z$,
$\ldots$ in  $D$, we always  denote by $x'$, $y'$, $z'$, $\ldots$
the images in $D'$ of $x$, $y$, $z$, $\ldots$ under $f$,
respectively. Also we assume that $\alpha$, $\beta$, $\gamma$,
$\ldots$ denote curves in $D$ and $\alpha'$, $\beta'$, $\gamma'$,
$\ldots$  the images in $D'$ of $\alpha$, $\beta$, $\gamma$,
$\ldots$ under $f$, respectively.

\subsection{Uniform domains}

In their paper \cite{MS},  Martio and Sarvas  introduced uniform
domains. Now there are many alternative characterizations for
uniform domains, see \cite{FW,  Geo, Martio-80, Vai6, Vai4, Vai8}.
The importance of this class of domains in the function theory is
well documented; see \cite{FW, Vai2} etc. Moreover, uniform domains
in $\mathbb{R}^n$ enjoy numerous geometric and function theoretic
features in many areas of modern mathematical analysis; see
\cite{Alv, Bea, Geo, Has, Ml, Jo80, Yli, Vai2}. We adopt the
definition of uniform domains following closely the notation and
terminology of \cite{TV, Vai2, Vai, Vai6-0, Vai6} or
\cite{Martio-80}.

\bdefe \label{def1.3} A domain $D$ in $E$ is called $c$-{\it
uniform} in the norm metric provided there exists a constant $c$
with the property that each pair of points $z_{1},z_{2}$ in $D$ can
be joined with a rectifiable arc $\alpha$ in $ D$ satisfying

 \bee
\item\label{wx-4} $\ds\min_{j=1,2}\ell (\alpha [z_j, z])\leq c\,d_{D}(z)
$ for all $z\in \alpha$, and

\item\label{wx-5} $\ell(\alpha)\leq c\,|z_{1}-z_{2}|$,
\eee

\noindent where $\ell(\alpha)$ denotes the length of $\alpha$,
$\alpha[z_{j},z]$ the part of $\alpha$ between $z_{j}$ and $z$, and
$d_D(z)$ the distance from $z$ to the boundary $\partial D$ of $D$.
At this time, $\alpha$ is said to be a {\it double $c$-cone arc}.
\edefe

\subsection{Quasihyperbolic metric, quasihyperbolic geodesics, solid arcs and neargeodesics}

Gehring and Palka \cite{GP} introduced the quasihyperbolic metric of
a domain in $\IR^n$. After that, this metric has become an important
tool in geometric function theory and in its generalizations to
metric spaces and to Banach spaces, see
 \cite{Avv, Bea, BHX, GP, Geo, HIMPS, KN, Vai6-0, Vai6, Vai8, Mvo1} etc. Yet, some basic questions
of the quasihyperbolic geometry in Banach spaces are open. For
instance, only recently the convexity of quasihyperbolic balls has
been studied in \cite{k, RT, Vai3} and the smoothness of
 geodesics in \cite{RT2} in the setup of Banach spaces.

For each pair of points $z_1$, $z_2$ in $D$, the {\it distance ratio
metric} $j_D(z_1,z_2)$ between $z_1$ and $z_2$ is defined by \cite{GP,vu85}
$$j_D(z_1,z_2)=\log\Big(1+\frac{|z_1-z_2|}{\min\{d_D(z_1),d_D(z_2)\}}\Big).$$

The {\it quasihyperbolic length} of a rectifiable arc or a path
$\alpha$ in the norm metric in $D$ is the number (cf.
\cite{GP,Vai3}):

$$\ell_{k_D}(\alpha)=\int_{\alpha}\frac{|dz|}{d_{D}(z)}.
$$

For each pair of points $z_1$, $z_2$ in $D$, the {\it quasihyperbolic distance}
$k_D(z_1,z_2)$ between $z_1$ and $z_2$ is defined in the usual way:
$$k_D(z_1,z_2)=\inf\ell_{k_D}(\alpha),
$$
where the infimum is taken over all rectifiable arcs $\alpha$
joining $z_1$ to $z_2$ in $D$. For all $z_1$, $z_2$ in $D$, we have
(cf. \cite{Vai3})

\beq\label{eq(0000)} k_{D}(z_1, z_2)&\geq&
\inf\left\{\log\Big(1+\frac{\ell(\alpha)}{\min\{d_{D}(z_1), d_{D}(z_2)\}}\Big)\right\}\\ \nonumber
&\geq& j_D(z_1,z_2)
\geq \Big|\log \frac{d_{D}(z_2)}{d_{D}(z_1)}\Big|,\eeq
where the infimum is taken over all rectifiable curves $\alpha$ in $D$ connecting $z_1$ and $z_2$.   Moreover, if $|z_1-z_2|\le d_D(z_1)$, we have
\cite{Vai6-0, vu85}
\begin{equation} \label{vu1}
k_D(z_1,z_2)\le \log\Big( 1+ \frac{
|z_1-z_2|}{d_D(z_1)-|z_1-z_2|}\Big)\le \frac{
|z_1-z_2|}{d_D(z_1)-|z_1-z_2|}\,,
\end{equation}
where the last inequality follows from the following elementary
inequality
\begin{equation} \label{logINE}
\frac{r}{1-r/2} \le \log \frac{1}{1-r} \le \frac{r}{1-r} \, \quad {\rm  for }\,\, 0\le r<1 \,.
\end{equation}

In \cite{Vai6},  V\"ais\"al\"a characterized uniform domains by the
quasihyperbolic metric.

\begin{Thm}\label{thm0.1} {\rm (\cite[Theorem 6.16]{Vai6})}
For a domain $D\not= E$, the following are quantitatively equivalent: \bee

\item $D$ is a $c$-uniform domain;
\item $k_D(z_1,z_2)\leq c'\;j_D(z_1,z_2)$ for all $z_1,z_2\in D$;
\item $k_D(z_1,z_2)\leq c'_1\;
 j_D(z_1,z_2)+d$ for all $z_1,z_2\in D$.\eee
 where $c$ and $c'$ depend on each other, $c$ and $c'_1$ depend on each other, and $d$ depends on $c$ (or $c'$).\end{Thm}

 In the case of domains in $ {\mathbb R}^n \,,$ the equivalence
  of items (1) and (3) in Theorem \Ref{thm0.1} is due to Gehring and Osgood \cite{Geo} and the
  equivalence of items (2) and (3) due to Vuorinen \cite[2.50 (2)]{vu85}.

Recall that an arc $\alpha$ from $z_1$ to
$z_2$ is a {\it quasihyperbolic geodesic} if
$\ell_{k_D}(\alpha)=k_D(z_1,z_2)$. Each subarc of a quasihyperbolic
geodesic is obviously a quasihyperbolic geodesic. It is known that a
quasihyperbolic geodesic between every pair of points in $E$ exists if the
dimension of $E$ is finite, see \cite[Lemma 1]{Geo}. This is not
true in arbitrary spaces (cf. \cite[Example 2.9]{Vai6-0}).
In order to remedy this shortage, V\"ais\"al\"a introduced the following concepts \cite{Vai6}.

\bdefe \label{def1.4}
 Let $\alpha$ be an arc in $E$. The arc
may be closed, open or half open. Let $\overline{x}=(x_0,...,x_n)$,
$n\geq 1$, be a finite sequence of successive points of $\alpha$.
For $h\geq 0$, we say that $\overline{x}$ is {\it $h$-coarse} if
$k_D(x_{j-1}, x_j)\geq h$ for all $1\leq j\leq n$. Let $\Phi_k(\alpha,h)$
be the family of all $h$-coarse sequences of $\alpha$. Set

$$s_k(\overline{x})=\sum^{n}_{j=1}k_D(x_{j-1}, x_j)$$ and
$$\ell_{k_D}(\alpha, h)=\sup \{s_k(\overline{x}): \overline{x}\in \Phi_k(\alpha,h)\}$$
with the agreement that $\ell_{k_D}(\alpha, h)=0$ if
$\Phi_k(\alpha,h)=\emptyset$. Then the number $\ell_{k_D}(\alpha, h)$ is the
{\it $h$-coarse quasihyperbolic length} of $\alpha$.\edefe

\bdefe \label{def1.5} Let $D$ be a domain in $E$. An arc $\alpha\subset D$
is {\it $(\nu, h)$-solid} with $\nu\geq 1$ and $h\geq 0$ if
$$\ell_{k_D}(\alpha[x,y], h)\leq \nu\;k_D(x,y)$$ for all $x, y\in \alpha$.
A {\it $(\nu,0)$-solid arc} is said to be a {\it $\nu$-neargeodesic}, i.e.
an arc $\alpha\subset D$ is a $\nu$-neargeodesic if and only if $\ell_{k_D}(\alpha[x,y])\leq \nu\;k_D(x,y)$
for all $x, y\in \alpha$.\edefe

Obviously, a $\nu$-neargeodesic is a quasihyperbolic geodesic if and
only if $\nu=1$.

In \cite{Vai6}, V\"ais\"al\"a got the following property concerning
the existence of neargeodesics in $E$.

\begin{Thm}\label{LemA} $($\cite[Theorem 3.3]{Vai6}$)$
Fix  $\nu>1\,.$  Then for
 $\{z_1,\, z_2\}\subset D \subset E$  there is a
$\nu$-neargeodesic in $D$ joining $z_1$ and $z_2$.
\end{Thm}


\begin{Thm}\label{Thm4-1} {\rm (\cite[Theorem 6.22]{Vai6})} Suppose that
$\gamma\subset D\not= E$ is a $(\nu, h)$-solid arc with endpoints $z_1$, $z_2$,
and that $D$ is a $c$-uniform domain. Then there is a constant $\mu_1=\mu_1(\nu,h,c)\geq 1$
such that
\bee
\item  $\ds\min\{\diam(\gamma [z_1, z]),\diam(\gamma [z_2, z])\}\leq \mu_1\,d_{G}(z)
$ for all $z\in \gamma$, and

\item  $\diam(\gamma)\leq \mu_1\max\big\{|z_1-z_2|,2(e^h-1)\min\{d_{G}(z_1),d_{G}(z_2)\}\big\}$.
\eee\end{Thm}

\subsection{FQC, QH and CQH mappings}

\bdefe \label{def1.7} Let $G\not=E$ and $G'\not=E'$ be metric
spaces, and let $\varphi:[0,\infty)\to [0,\infty)$ be a growth
function, that is, a homeomorphism with $\varphi(t)\geq t$. We say
that a homeomorphism $f: G\to G'$ is {\it $\varphi$-semisolid} if
$$ k_{G'}(f(x),f(y))\leq \varphi(k_{G}(x,y))$$
for all $x$, $y\in G$, and {\it $\varphi$-solid} if both $f$ and $f^{-1}$
satisfy this condition.

We say that $f$ is {\it fully $\varphi$-semisolid}
(resp. {\it fully $\varphi$-solid}) if $f$ is
$\varphi$-semisolid (resp. $\varphi$-solid) on every  subdomain of $G$. In particular,
when $G=E$, the corresponding subdomains are taken to be proper ones. Fully $\varphi$-solid mappings are also called {\it freely
$\varphi$-quasiconformal mappings}, or briefly {\it $\varphi$-FQC mappings}.\edefe

If $E=\IR^n=E'$, then $f$ is $FQC$ if and only if $f$ is
quasiconformal (cf. \cite{Vai6-0}). See \cite{Vai1, Mvo1} for definitions and
properties of $K$-quasiconformal mappings, or briefly $K$-QC mappings.

\bdefe \label{def1.6} We say that a homeomorphism $f: D\rightarrow
D'$, where $D\subset E$ and $D'\subset E'$, is {\it $C$-coarsely $M$-quasihyperbolic}, or briefly
$(M,C)$-CQH, in the quasihyperbolic metric if it satisfies
$$\frac{k_D(x,y)-C}{M}\leq k_{D'}(f(x),f(y))\leq M\;k_D(x,y)+C$$
for all $x$, $y\in D$. An $(M,0)$-CQH mapping is $M$-bilipschitz in the quasihyperbolic metric,
or briefly $ M$-QH. \edefe

The following result easily follows from the definitions.

\begin{prop}\label{fin-1}
Every $\varphi$-FQC mapping is an $(M,C)$-CQH mapping, where $M$ and $C$
depend only on $\varphi$.\end{prop}

\begin{Thm}\label{Lem4-0} {\rm (\cite[Theorem 4.15]{Vai6})} For domains
$D\not= E$ and $D'\not=E'$, suppose that $f: D\to D'$ is $(M,C)$-CQH.
If $\gamma$ is a $(\nu,h)$-solid arc in $D$, then the image arc
$\gamma'$ is $(\nu',h_1)$-solid in $D'$ with $(\nu', h_1)$
depending only on $(M, C, \nu,h)$.\end{Thm}

\medskip

\subsection{Quasisymmetric mappings}
Let $X$ be a metric space and $\dot{X}=X\cup \{\infty\}$. By a triple in $X$ we mean an ordered sequence
$T=(x,a,b)$ of three distinct points in $X$. The ratio of $T$ is the number $$\rho(T)=\frac{|a-x|}{|b-x|}.$$
If $f: X\to Y$ is  an injective map, the image of a triple  $T=(x,a,b)$  is the triple $fT=(fx,fa,fb)$.

Suppose that  $A\subset X$. A triple  $T=(x,a,b)$ in $X$ is said to be a triple
in the pair $(X, A)$ if $x\in A$ or if $\{a,b\}\subset A$. Equivalently, both $|a-x|$ and $|b-x|$ are distances from a point in $A$.

\bdefe \label{def1-0} Let $X$ and $Y$ be two metric spaces,
and let $\eta: [0, \infty)\to [0, \infty)$ be a homeomorphism. Suppose
$A\subset X$. An embedding $f: X\to Y$ is said to be {\it
$\eta$-quasisymmetric} relative to $A$, or briefly $\eta$-$QS$ rel
$A$, if $\rho(T)\leq t$ implies that $\rho(f(T))\leq \eta(t)$  for each
triple $T$ in $(X,A)$ and $t\geq 0$. \edefe
Thus ``quasisymmetry rel $X$" is equivalent to ordinary ``quasisymmetry".

\bdefe \label{def1-0'} Let  $0<q<1$, let $\eta$ be as in Def. \ref{def1-0} and let $G$, $G'$ be metric
spaces in $E$ and $E'$, respectively. A homeomorphism $f:G\to G'$ is
{\it $q$-locally $\eta$-quasisymmetric} if $f|_{\mathbb{B}(a,qr)}$ is
$\eta$-QS whenever $\mathbb{B}(a,r)\subset G$. If $G\not=E$, this
means that $f|_{\mathbb{B}(a, q d_G(a))}$ is $\eta$-QS. When $G=E$, obviously, $f$ is $\eta$-QS. \edefe

It is known that each $K$-QC  mapping in
$\IR^n$ is $q$-locally $\eta$-QS for every $q<1$ with $\eta=\eta(K,
q, n)$, i.e. $\eta$ depends only on the constants $K$, $q$ and $n$
(cf. \cite[5.23]{Avv}). Conversely, each $q$-locally $\eta$-QS
mapping in $\IR^n$ is a $K$-QC mapping with $K=(\eta(1))^{n-1}$ by
the metric definition of quasiconformality (cf. \cite[5.6]{Vai6-0}).
Further, in \cite{Vai6-0}, V\"ais\"al\"a proved

\begin{Thm}\label{Thm1}$($\cite[Theorem 5.10]{Vai6-0}$)$ Suppose that
$D$ and
$D'$ are domains in $E$ and $E'$, respectively. For a homeomorphism $f: D\to D'$, the following conditions are quantitatively equivalent:

\bee
\item $f$ is $\psi$-FQC;
\item both $f$ and $f^{-1}$ are  $q$-locally $\eta$-QS;
\item For every $0<q<1$, there is some  $\eta(q)$ such that both $f$ and $f^{-1}$ are  $q$-locally $\eta(q)$-QS.
\eee\end{Thm}

\begin{Thm}\label{Thm2'}$($\cite[Theorem 5.13]{Vai6-0}$)$
Suppose that $f: E\to D'\subset E'$ is fully $\varphi$-semisolid. Then the following conditions are quantitatively equivalent:
\bee
\item $D'=E'$;
\item $f$ is $\eta$-QS with $\eta=\eta(\varphi)$;
\item $f$ is $\psi$-FQC with $\psi=\psi(\varphi)$.
\eee\end{Thm}

\section{The proof of Theorem \ref{thm1.1}}\label{sec-4}


Before the proof of Theorem \ref{thm1.1} we need some preparation.\medskip

\noindent{\bf Basic assumption}\quad In the following, we always
assume that $f: D\to D'$ is $\psi$-FQC and that $D'$ is an $a$-uniform
domain. By Proposition \ref{fin-1}, we assume further that both $f$ and $f^{-1}$ are $(M, C)$-CQH homeomorphisms,
where
$M$ and $C$ depend only on $\psi$.
We see from Theorem \Ref{Thm2'} that
$D=E$ if and only if $D'=E'$. Hence, to prove Theorem \ref{thm1.1}, it suffices to consider the case
$D\not=E$ and $D'\not=E'$.  Since $D'$ is $a$-uniform, Theorem \Ref{thm0.1} implies that there is a constant $a'$ such that

\beq\label{Wen-1} k_{D'}(\xi', \zeta')\leq a'\;
 j_{D'}(\xi', \zeta')\eeq for every pair
 of points $\xi'$, $\zeta'\in D'$, where
 $a'\leq 7a^3$ (cf. \cite[Theorem 2.23]{Vai4}).
 Theorem \Ref{Thm1} shows
that both $f$ and $f^{-1}$ are $\frac{3}{4}$-locally $\eta$-QS with
$\eta$ depending on $\psi$. Without loss of generality, we may assume that $\eta(1)\geq 1$. Let $D_1\subset D$ be a
$c$-uniform domain. Then it follows from Theorem \Ref{thm0.1} that there is a constant $c'$ such that

\beq\label{Wen-2} k_{D_1}(u, v)\leq c'\;j_{D_1}(u, v)\eeq
for every pair of points $u,\, v\in D_1\,,$ where
 $c'\leq 7c^3$ (cf. \cite[Theorem 2.23]{Vai4}).
Since $D'$ is $a$-uniform, without
loss of generality, we may assume that $D_1$ is a proper subdomain
of $D$.
 For a pair of points $x'$, $y'\in D'_1$,
let $\gamma'$ be a $2$-neargeodesic  in $D'_1$ joining $x'$ and
$y'$. It follows from Theorem \Ref{Lem4-0} that $\gamma$ is
a $(\nu',h_1)$-solid arc joining $x$ and $y$ in $D_1$, where $(\nu',h_1)$ depends only on $(M,C)$.
\medskip

Let $z_0\in\gamma$ be such that

$$d_{D_1}(z_0)=\sup\limits_{p\in \gamma}d_{D_1}(p).
$$
It is possible that $z_0=x$ or $y$. Then

\begin{Thm}\label{lem-j-j} $($\cite[Lemma 2.1]{HLVW}$)$ There is a constant ${\rho}_1=
\rho_1(c,\nu', h_1)$ with $\rho_1\geq 1$ such that

\noindent $(1)$ For all  $z\in
\gamma[x, z_0]$,
 $$|x-z|\leq
\rho_1\;d_{D_1}(z),$$ and for each  $z\in \gamma[y, z_0]$,
$$|y-z|\leq \rho_1\;d_{D_1}(z);$$
\noindent $(2)$ $\diam(\gamma)\leq \rho_1\max\big\{|x-y|,2(e^{h_1}-1)\min\{d_{D_1}(x),d_{D_1}(y)\}\big\}.$
\end{Thm}

By Theorem $H$, we get

\bcor\label{cor1}
\noindent $(1)$ For each $x_1\in
\gamma[x, z_0]$ and for all $z\in \gamma[x_1, z_0]$,
 $$|x_1-z|\leq
\rho_1\;d_{D_1}(z);$$

\noindent $(2)$ For each $y_1\in
\gamma[y, z_0]$ and for all $z\in \gamma[y_1, z_0]$,
 $$|y_1-z|\leq
\rho_1\;d_{D_1}(z),$$

\noindent where $\rho_1$ is the same as in Theorem $H$.
\ecor

For the convenience of the statements and proofs of the lemmas below, we write down the related constants:
\medskip

\noindent $(1)$ $\vartheta=\min\{\frac{1}{2}\eta^{-1}(\frac{1}{\rho_1}),
\frac{1}{2\rho_1}\}$,

\noindent $(2)$
$b_1=\max\big\{\frac{1}{\psi^{-1}(\frac{1}{8})} , 1 \big\} 4^{8a'c'CM\rho_1\eta(\rho_1)\psi(1)}$,

\noindent $(3)$ $b_2= \frac{1}{\eta^{-1}(\vartheta)}b_1^{8a'c'CM}$,

\noindent $(4)$ $b_3= 2e^{5a'c'CMb_2^2}$,

\noindent $(5)$
$b_4=\frac{e^{\vartheta_1}(\eta(13b_3^2))^{3a'M}}{2a'c'CM\eta^{-1}(\vartheta)}$,

\noindent where $a'$, $c'$, $C$, $M$, $\eta$ and $\psi$ are the same
as in {\bf Basic assumption}, $\rho_1$ is the constant in Theorem $H$ and
$\vartheta_1=\frac{5a'b_3^4c'M^2}{\eta^{-1}(\frac{1}{b_3^3})}$.
Obviously, $b_k>1$ for each $k\in\{1,2,3,4\}$.
\medskip

Since $\eta(1)\geq 1$, we also easily get
\begin{cor}\label{dear-1} $(1)$ $b_2>4^{(8a'c'CM)^2\rho_1}$;

$(2)$ $b_3>(4a'b_2c'M^2)^{6a'b_2c'CM^2+2}$;

$(3)$ $b_4>\max\{24b_3^2, (7b_3^2)^{8a'c'M^2\tau+2}\}$, where $\tau=\frac{2}{\eta^{-1}(\frac{1}{b_3^3})}$. \end{cor}

 In what follows, we prove that $\gamma'$ is a $b_4$-double cone arc in $D'_1$. (Recall that $\gamma'$
 is a $2$-neargeodesic in $D_1'$ joining $x'$ and $y'$ as constructed in Basic
 assumption.)
 More precisely, we will prove
that $\gamma'$ satisfies the following conditions:

\be\label{lem-j-j}\min\{\ell(\gamma'[x',z']),\ell(\gamma'[y',z'])\}\leq
2b_3^2\;d_{D'_1}(z')\,\; \mbox{for}\,\mbox{all}\,\; z'\in \gamma'\ee and

\be\label{lem0-1} \ell(\gamma')\leq b_4|x'-y'|.\ee

In the remaining part of this paper our main task will be to prove (\ref{lem-j-j})
and (\ref{lem0-1}).

\subsection{ The proof of \eqref{lem-j-j}.}

If for each $z'\in
\gamma'[x',z'_0]$,
$$\ell(\gamma'[x',z'])\leq b_3^2\,d_{D'_1}(z'),$$ and for
 each $z'\in
\gamma'[y',z'_0]$,
$$\ell(\gamma'[y',z'])\leq b_3^2\,d_{D'_1}(z'),$$ then \eqref{lem-j-j} is obvious. Hence
the remaining cases we need to consider are: either there is some $p'\in
\gamma'[x',z'_0]$ such that $\ell(\gamma'[x',p'])> b_3^2\,d_{D'_1}(p')$ or there is some $q'\in
\gamma'[y',z'_0]$ such that $\ell(\gamma'[y',q'])> b_3^2\,d_{D'_1}(q')$. We only need to consider the former case,
that is, there is some $p'\in
\gamma'[x',z'_0]$ such that
\be\label{hl-1} \ell(\gamma'[x',p'])> b_3^2\,d_{D'_1}(p'),\ee
since the argument for the latter one is similar.

We use $w'_0$ to denote the first point in $\gamma'[x',z'_0]$ in the direction
from $x'$ to $z'_0$ such that
\be\label{eq(3-1)}\ell(\gamma'[x',w'_0])=b_3\,d_{D'_1}(w'_0).\ee
Inequality \eqref{hl-1} guarantees the existence of such a point $w_0'$.

\begin{lem} \label{lem-3-1-1}  For every $z\in\gamma[w_0,z_0]$, we have
$d_{D}(z)\leq b_2d_{D_1}(z)$.
\end{lem}
\bpf We prove this lemma by contradiction. Suppose that there exists
some point $w_1\in\gamma[w_0,z_0]$ such that

\be\label{hw-1} d_{D}(w_1)> b_2d_{D_1}(w_1).\ee
Then it follows from Corollary \ref{cor1}
that for all $w\in \gamma[x,w_1]$,

\be\label{hu-3-1}|w-w_1|\leq
 \rho_1\,d_{D_1}(w_1)<\frac{\rho_1}{b_2}d_{D}(w_1),\ee
and hence by \eqref{vu1} and Corollary \ref{dear-1}
$$k_{D}(w,w_1) \leq
\frac{\rho_1}{b_2-\rho_1},$$
and so
$$k_{D'}(w',w'_1)\leq \psi(k_{D}(w,w_1))\leq \psi\Big(\frac{\rho_1}{b_2-\rho_1}\Big)\leq \frac{1}{8},$$ since $f$ is $\psi$-FQC. We have

\be\label{hm-2}|w'-w'_1|\leq \frac{1}{7}d_{D'}(w'_1),\ee because otherwise by the lower bound of \eqref{logINE},
$k_{D'}(w',w'_1)\geq j_{D'}(w',w'_1)\ge \frac{2}{15}\ge \frac{1}{8}\,.$

The triangle inequality yields
$d_{D'}(w'_0)\geq d_{D'}(w'_1)-|w'_0-w'_1|$ and further by \eqref{hm-2} we see that
$$d_{D'}(w'_0)\geq \frac{6}{7}d_{D'}(w'_1).$$ Then for  all $w'\in \gamma'[x',w_1']$, we have

\be\label{xl-1}|w'-w'_0|\leq |w'-w'_1|+|w'_1-w'_0|\leq \frac{2}{7}d_{D'}(w'_1)\leq \frac{1}{3}d_{D'}(w'_0).\ee

Let $w'_2\in \partial D'_1$ be such that $|w'_2-w'_0|\leq d_{D'_1}(w'_0)+\frac{1}{6}\min\{d_{D'}(w'_0),d_{D'_1}(w'_0)\}$,
and let $w'_3\in \gamma'[x',w'_0]$ be such that
$\ell(\gamma'[x',w'_3])= \frac{1}{2}\ell(\gamma'[x',w'_0])$. It
follows from the assumption ``$\gamma'$ being $2$-neargeodesic",
\eqref{eq(0000)} and \eqref{eq(3-1)} that

\beq\label{hw-2} k_{D'_1}(w'_0,w'_3)
&\geq & \frac{1}{2}\ell_{k_{D'_1}}(\gamma'[w'_3,w'_0])\\ \nonumber
&\geq& \frac{1}{2}\log\left(1+
\frac{\ell(\gamma'[w_3', w'_0])}{\min\{d_{D_1'}(w_3'), d_{D_1'}(w_0')\}}\right)\\ \nonumber
&\geq& \frac{1}{2}\log\Big(1+\frac{b_3}{2}\Big).\eeq

Suppose that $d_{D'_1}(w'_0)\geq \frac{3}{2}|w'_0-w'_3|.$ Then by
\eqref{vu1}
$k_{D'_1}(w'_0,w'_3)\leq 2$,
 which contradicts \eqref{hw-2}.
Hence we have proved that

\be\label{h-l-5} d_{D'_1}(w'_0)< \frac{3}{2}|w'_0-w'_3|.\ee
We infer from
(\ref{xl-1}) and \eqref{h-l-5} that

\beq\label{hm-3}|w'_2-w'_0|<\frac{3}{2}|w'_0-w'_3|+\frac{1}{6}d_{D'}(w'_0)<\frac{2}{3}d_{D'}(w'_0).\eeq
\begin{figure}[!ht]
\centering
\includegraphics[width=9cm]{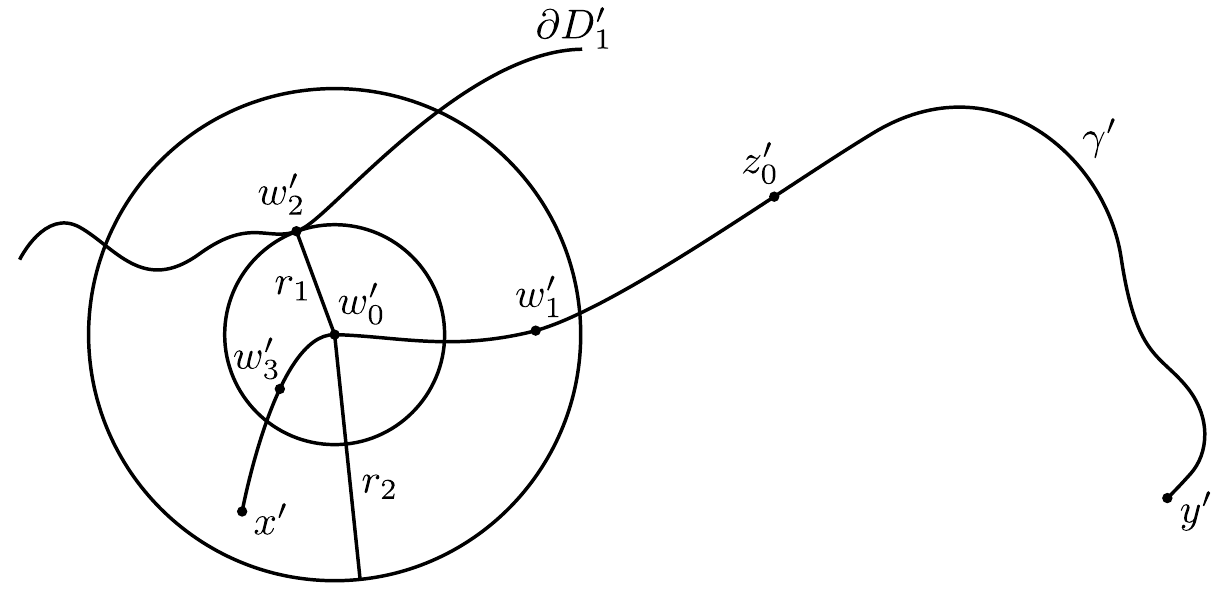}
\caption{\label{fig4}  $r_1\leq \frac{7}{6}d_{D'_1}(w'_0), r_2=\frac{2}{3}d_{D'}(w'_0)$}
\end{figure}
Thus (\ref{xl-1}) and (\ref{hm-3}) imply that $w'_2, w'_3\in \mathbb{B}(w'_0,\frac{2}{3}d_{D'}(w'_0))$.
We see from Corollary \ref{cor1}, the choice of $w'_2$ and the assumption ``$f^{-1}$ being $\frac{3}{4}$-locally
$\eta$-QS" that
$$\frac{1}{\rho_1}\leq\frac{|w_0-w_2|}{|w_0-w_3|}\leq \eta\Big(\frac{|w'_0-w'_2|}{|w'_0-w'_3|}\Big).$$
Then the choice of $w_0'$ and $w_3'$ yields that
\be\label{xl-2}|w'_0-w'_3|\leq \frac{1}{\eta^{-1}(\frac{1}{\rho_1})}|w'_0-w'_2|\leq \frac{7}{6\eta^{-1}(\frac{1}{\rho_1})}d_{D'_1}(w'_0)\leq \frac{7}{3\eta^{-1}(\frac{1}{\rho_1})}d_{D'_1}(w'_3).\ee
Moreover, we infer from \eqref{Wen-2} and
\eqref{hw-2} that

\begin{eqnarray*}j_{D_1}(w_0,w_3)&\geq& \frac{1}{c'}
k_{D_1}(w_0,w_3) \geq \frac{1}{c'M}\big(k_{D'_1}(w'_0,w'_3)-C\big)
\\ \nonumber &>& \log\Big(1+ e^{a'b_2C}\Big).\end{eqnarray*}
Hence

\be\label{eq(3-2)}|w_0-w_3|\geq e^{a'b_2C}\min\{d_{D_1}(w_3),d_{D_1}(w_0)\}> b_2\,d_{D_1}(w_3),\ee
since by Corollary \ref{cor1}, $d_{D_1}(w_3)\leq d_{D_1}(w_0)+|w_0-w_3|\leq (1+\rho_1)d_{D_1}(w_0)$.

\begin{figure}[!ht]
\centering
\includegraphics[width=8cm]{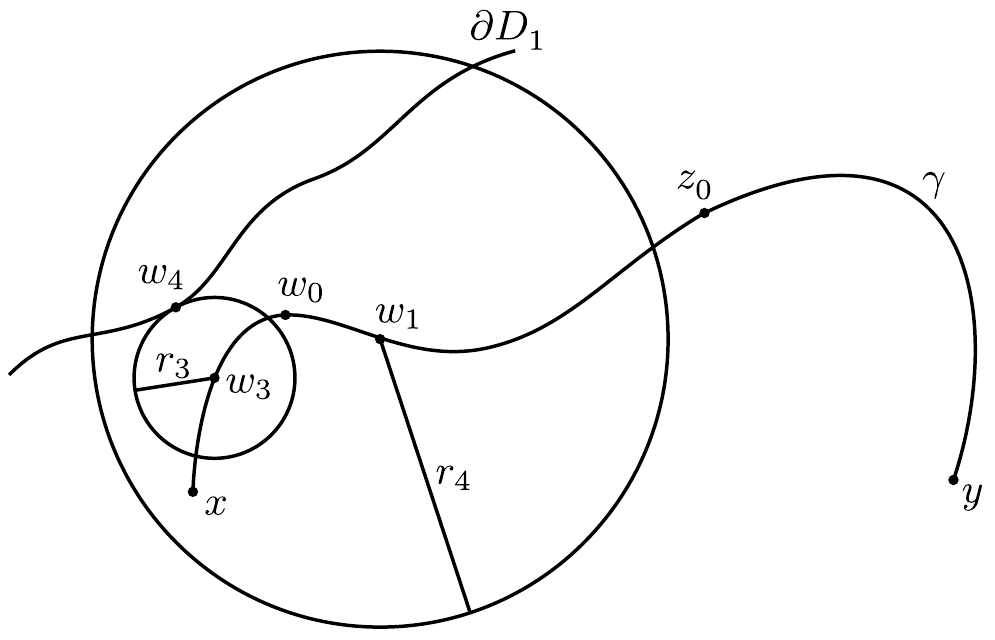}
\caption{\label{fig3} $ r_3\leq 2d_{D_1}(w_3), r_4=\frac{1}{2}d_{D}(w_1)$}
\end{figure}
Let $w_4\in\partial D_1$ be such that $|w_4-w_3|\leq d_{D_1}(w_3)+\min\{|w_3-w_1|,d_{D_1}(w_3)\}$. Then
we obtain from \eqref{hw-1} and (\ref{hu-3-1}) that

\beq\label{hm-5}|w_4-w_1|&\leq&
|w_4-w_3|+|w_3-w_1| < 3|w_3-w_1|+d_{D_1}(w_1)\\
\nonumber&\leq&\frac{1+3\rho_1}{b_2}d_{D}(w_1),\eeq
since $|w_4-w_3|\leq d_{D_1}(w_3)+|w_3-w_1|< d_{D_1}(w_1)+2|w_3-w_1|$. Hence (\ref{hu-3-1}) and (\ref{hm-5}) show that $w_0, w_3, w_4\in \mathbb{B}(w_1,\frac{1}{2}d_{D}(w_1))$.
It follows from \eqref{xl-2} and (\ref{eq(3-2)}) that
$$\frac{3\eta^{-1}(\frac{1}{\rho_1})}{7}\leq\frac{|w'_4-w'_3|}{|w'_0-w'_3|}\leq
\eta\Big(\frac{|w_4-w_3|}{|w_0-w_3|}\Big)\leq
\eta\Big(\frac{2}{b_2}\Big)<\frac{3\eta^{-1}(\frac{1}{\rho_1})}{7},
$$ since $f$ is $\frac{3}{4}$-locally $\eta$-QS.
This obvious contradiction completes the proof of Lemma \ref{lem-3-1-1}.
\epf

If $\ell(\gamma'[y',z'])\leq b_3\;d_{D'_1}(z')$ for all $z'\in
\gamma'[y',z'_0]$, then we let $y'_0=z'_0$. Otherwise, we let
$y'_0$ be the first point in $\gamma'[y',z'_0]$ in the direction
from $y'$ to $z'_0$ such that
\be\label{eq(3-4)}\ell(\gamma'[y',y'_0])= b_3\,d_{D'_1}(y'_0).\ee

A reasoning  similar to the one in the proof of Lemma \ref{lem-3-1-1} shows that

\begin{lem} \label{lem-3-1-1'}For all $z\in\gamma[y_0,z_0]$, we have  $d_{D}(z)\leq b_2d_{D_1}(z).$\end{lem}

Let $u'_0\in\gamma'[w'_0,y'_0]$ satisfy
$$d_{D'_1}(u'_0)=\sup_{z'\in\gamma'[w'_0,y'_0]}d_{D'_1}(z').$$
Obviously, there exists a nonnegative integer $m_1$
such that

$$ 2^{m_1}\, d_{D'_1}(w'_0) \leq d_{D'_1}(u'_0)< 2^{m_1+1}\, d_{D'_1}(w'_0).$$
Let $v'_0$ be the first point in $\gamma'[w'_0,u'_0]$ from $w'_0$ to
$u'_0$ with

$$d_{D'_1}(v'_0)=2^{m_1}\, d_{D'_1}(w'_0),
$$ and let $u_1'=w'_0$. If $v'_0= u'_1$, we let $u'_2=u'_0$. It is possible
that $u'_1=u'_2$. If $v'_0\not= u'_1$, then we let $u'_2,\ldots
,u'_{m_1+1}\in \gamma'[w'_0,v'_0]$ be the points such that for each
$i\in \{2,\ldots,m_1+1\}$, $u'_i$ denotes the first point from
$w'_0$ to $u'_0$ with
$$d_{D'_1}(u'_i)=2^{i-1}\, d_{D'_1}(u'_1).
$$
Obviously, $u'_{m_1+1}=v'_0$. If $v'_0\not= u'_0$, we use
$u'_{m_1+2}$ to denote $u'_0$. Then we have the following assertion.

\begin{lem}\label{eq(3-6)} For all
$i\in\{1,\cdots, m_1+1\}$ and $z'\in
\gamma'[u'_i,u'_{i+1}]$, we have
$$\ell(\gamma'[u'_i,u'_{i+1}])\leq
b_3\,d_{D'_1}(z').$$\end{lem}

\bpf
Lemma \ref{lem-3-1-1}, Lemma \ref{lem-3-1-1'}, \eqref{Wen-1}, \eqref{Wen-2} and the assumption ``
both $f$ and $f^{-1}$ being $(M, C)$-CQH" imply that
\begin{eqnarray*}\frac{\ell(\gamma'[u'_i,u'_{i+1}])}{2d_{D'_1}(u'_i)}&\leq&
\ell_{k_{D'_1}}(\gamma'[u'_i,u'_{i+1}])\leq 2k_{D'_1}(u'_i,u'_{i+1})
\leq  2Mk_{D_1}(u_i,u_{i+1})+2C
\\ \nonumber&\leq& 2c'Mj_{D_1}(u_i,u_{i+1})+2C
\\ \nonumber&\leq& 2c'M\log\Big(1+\frac{b_2|u_i-u_{i+1}|}{\min\{d_{D}(u_i),d_{D}(u_{i+1})\}}\Big)+2C
\\ \nonumber&\leq& 2b_2c'Mk_{D}(u_i,u_{i+1})+2C
\\ \nonumber&\leq& 2b_2c'M^2k_{D'}(u'_i,u'_{i+1})+2b_2c'CM+2C
\\ \nonumber&\leq& 2a'b_2c'M^2\log\Big(1+\frac{|u'_i-u'_{i+1}|}{d_{D'_1}(u'_i)}\Big)+2b_2c'CM+2C
\\ \nonumber&\leq& 2a'b_2c'M^2\log\Big(1+\frac{\ell(\gamma'[u'_i,u'_{i+1}])}{d_{D'_1}(u'_i)}\Big)+2b_2c'CM+2C,\end{eqnarray*}
and hence we easily deduce the following conclusion:
\be\label{eq(3-7-1)}
\ell(\gamma'[u'_i, u'_{i+1}])\leq (3a'b_2c'M^2)^2\,d_{D'_1}(u'_i),\ee
and we also see that for all $z'\in \gamma'[u'_i,u'_{i+1}]$,

\begin{eqnarray*}\log \frac{d_{D'_1}(u'_i)}{d_{D'_1}(z')}&\leq&
k_{D'_1}(u'_i,z')
\leq \ell_{k_{D'_1}}(\gamma'[u'_i,z'])
\leq 2k_{D'_1}(u'_i,u'_{i+1})
\\ \nonumber&\leq& 2Mk_{D_1}(u_i,u_{i+1})+2C
\\ \nonumber&\leq& 2c'M\;j_{D_1}(u_i,u_{i+1})+2C
\\ \nonumber&\leq& 2c'M\log\Big(1+\frac{b_2|u_i-u_{i+1}|}{\min\{d_{D}(u_i),d_{D}(u_{i+1})\}}\Big)+2C
\\ \nonumber&\leq& 2b_2c'Mk_{D}(u_i,u_{i+1})+2C
\\ \nonumber&\leq& 2b_2c'M^2k_{D'}(u'_i,u'_{i+1})+2b_2c'CM+2C
\\ \nonumber&\leq& 2a'b_2c'M^2\log\Big(1+\frac{|u'_i-u'_{i+1}|}{d_{D'_1}(u'_i)}\Big)+2b_2c'CM+2C
\\ \nonumber&\leq& 2a'b_2c'M^2\log\Big(1+\frac{\ell(\gamma'[u'_i,u'_{i+1}])}{d_{D'_1}(u'_i)}\Big)+2b_2c'CM+2C
\\ \nonumber&\leq& 2a'b_2c'M^2\log\Big(1+(3a'b_2c'M^2)^2\Big)+2b_2c'CM+2C.\end{eqnarray*} It follows that

\be\label{xt-1-1} d_{D'_1}(u'_i)\leq (4a'b_2c'M^2)^{ 6a'b_2c'CM^2}d_{D'_1}(z').\ee Hence the proof of Lemma
\ref{eq(3-6)} is complete by the combination of (\ref{eq(3-7-1)}), \eqref{xt-1-1} and Corollary \ref{dear-1}.\epf

\noindent {\bf Now we are ready to finish the proof of \eqref{lem-j-j}}.\medskip

For all $z'\in\gamma'[x',u'_0]$, if $z'\in\gamma'[x',w'_0]$, then \eqref{lem-j-j} easily follows from the choice of
$w_0'$.
For the case $z'\in\gamma'[w'_0, u'_0]$, we know that there exists some $k\in
\{1,\cdots,m_1+1\}$ such that $z'\in \gamma'[u'_k, u'_{k+1}]$. If $k=1$,
then by (\ref{eq(3-1)}), (\ref{xt-1-1}) and Lemma
\ref{eq(3-6)}, we know

\begin{eqnarray*}
\ell(\gamma'[x',z'])&=& \ell(\gamma'[x',w'_0])+\ell(\gamma'[w'_0,z'])\\
\nonumber &\leq& b_3(b_3+1)d_{D'_1}(z').\end{eqnarray*}

 If $k>1$,
again, we obtain from (\ref{eq(3-1)}), (\ref{xt-1-1}) and Lemma
\ref{eq(3-6)} that

\begin{eqnarray*}
\ell(\gamma'[x',z'])&=& \ell(\gamma'[x',w'_0])+\ell(\gamma'[w'_0,z'])\\
\nonumber &\leq& \ell(\gamma'[x',w'_0])+\ell(\gamma'[u'_1,u'_2])+\cdots+\ell(\gamma'[u'_k,u'_{k+1}])\\
\nonumber &\leq& b_3\big(2d_{D'_1}(u'_1)+d_{D'_1}(u'_2)+\cdots+d_{D'_1}(u'_k)\big)\\
\nonumber &\leq& 2{b_3}^2d_{{D'}_1} (z').\end{eqnarray*}
Hence we have proved that for all $z'\in\gamma'[x', u'_0]$,
$$ \ell(\gamma'[x',z'])\leq
2b_3^2\;d_{D'_1}(z').$$

Similar arguments as above show that
for all $z'\in \gamma'[y',u'_0]$,
$$\ell(\gamma'[y',z'])\leq 2b_3^2\,d_{D'_1}(z').$$
Hence the proof of \eqref{lem-j-j} is finished.

By \eqref{lem-j-j}, we have

\bcor\label{cor-2}
For every $x_1'\in \gamma'[x', y']$ and all $y_1'\in \gamma'[x'_1, y']$, we have that for each $z'\in \gamma'[x'_1,y'_1]$,
$$\min\{\ell(\gamma'[x'_1, z']), \ell(\gamma'[y_1', z'])\}\leq
2b_3^2\;d_{D'_1}(z').$$
\ecor

\subsection{The proof of \eqref{lem0-1}}

Suppose on the contrary that
\be\label{ly-0}\ell(\gamma')> b_4|x'-y'|.\ee  Then we claim $$ |x'-y'|>\frac{3}{4}\max\{d_{D'_1}(x'),d_{D'_1}(y')\}.$$

If
$$|x'-y'|\leq \frac{3}{4}\max\{d_{D'_1}(x'),d_{D'_1}(y')\},$$
then without loss of generality, we assume that
$$\max\{d_{D'_1}(x'),d_{D'_1}(y')\}=d_{D'_1}(x'). $$
Hence

$$d_{D'_1}(y')\geq
\frac{1}{4}d_{D'_1}(x'),$$
and by \eqref{eq(0000)} and \eqref{vu1}

\beq\label{my-2}
\log\Big(1+\frac{\ell(\gamma'[x',y'])}{d_{D'_1}(y')}\Big)&\leq& \ell_{k_{D'_1}}(\gamma'[x',y'])
\leq 2k_{D'_1}(x',y')
\\ \nonumber &\leq& \frac{8|x'-y'|}{d_{D'_1}(x')} \le 6 \,,
\eeq
since $d_{D'_1}(w')\geq d_{D'_1}(x')-|x'-w'|\geq \frac{1}{4}d_{D'_1}(x')$ for all $w'\in [x', y']$.
Thus \eqref{my-2} implies

\begin{eqnarray*} \frac{\ell(\gamma'[x',y'])}{e^6d_{D'_1}(y')}&\leq&
\log\Big(1+\frac{\ell(\gamma'[x',y'])}{d_{D'_1}(y')}\Big)\leq
\frac{8|x'-y'|}{d_{D'_1}(x')}.\end{eqnarray*}
It follows that
$$\ell(\gamma'[x',y'])\leq 8e^6|x'-y'|,$$ which contradicts \eqref{ly-0}.
Hence we get
\be\label{eq(3-7')} |x'-y'|>\frac{3}{4}\max\{d_{D'_1}(x'),d_{D'_1}(y')\}.\ee

Suppose (\ref{lem0-1}) does not hold, that is we have \eqref{ly-0}. Then  Corollary \ref{dear-1} shows that there exist $v'_1\in\gamma'$, $v'_2\in \gamma'[y',v'_1]$ such that

\be\label{xt-mz-1}\ell(\gamma'[x',v'_1])=12b_3^2|x'-y'|\,\; \mbox{and}\;
\,\ell(\gamma'[y',v'_2])=12b_3^2|x'-y'|.\ee Then we have

\begin{lem}\label{xt-mz-11}$|x'-v'_1|\geq\frac{1}{2}d_{D'_1}(v'_1)$.\end{lem}

\bpf Suppose on the contrary that
$$|x'-v'_1|<\frac{1}{2}d_{D'_1}(v'_1).$$ Then by \eqref{lem-j-j} and (\ref{xt-mz-1}), we get
 $$d_{D'_1}(x')\geq d_{D'_1}(v'_1)-|x'-v'_1|\geq \frac{1}{2}d_{D'_1}(v'_1)\geq \frac{1}{6b_3^2}\ell(\gamma'[x',v'_1])= 2|x'-y'|,$$ which contradicts
(\ref{eq(3-7')}). Hence the proof of Lemma \ref{xt-mz-11} is complete.\epf

Similarly, we have

\bcor\label{xt-mz-11'}$|y'-v'_2|\geq\frac{1}{2}d_{D'_1}(v'_2)$.\ecor

It follows from \eqref{lem-j-j}, (\ref{xt-mz-1}) and  Lemma \ref{xt-mz-11} that

\be\label{xt-mz-12} \min\{d_{D'_1}(v'_1),d_{D'_1}(v'_2)\}\leq d_{D'_1}(v'_1)\leq 2|x'-v'_1|\leq 24b_3^2|x'-y'|\ee and

\be\label{ly-1}|x'-v'_1|\geq\frac{1}{2}d_{D'_1}(v'_1)\geq\frac{1}{6b_3^2}\ell(\gamma'[x',v'_1])= 2|x'-y'|,\ee
whence \eqref{Wen-2}, (\ref{ly-0}), \eqref{xt-mz-1} and \eqref{xt-mz-12} show that

\begin{eqnarray*} c'\;j_{D_1}(v_1,v_2)
&\geq& k_{D_1}(v_1,v_2) \geq
\frac{1}{M}k_{D'_1}(v'_1,v'_2)-\frac{C}{M}
\\ \nonumber
&\geq& \frac{1}{2M}\ell_{k_{D_1'}}(\gamma'[v'_1,v'_2])-\frac{C}{M}\\ \nonumber
&\geq& \frac{1}{2M}\log\Big(1+\frac{\ell(\gamma'[v'_1,v'_2])}{\min\{d_{D'_1}(v'_1),d_{D'_1}(v'_2)\}}\Big)
-\frac{C}{M} \\ \nonumber
&\geq&
\frac{1}{2M}\log\Big(1+\frac{b_4-24b_3^2}{24b_3^2}\Big)-\frac{C}{M}
\\ \nonumber
&>&\frac{1}{3M}\log\Big(1+\frac{b_4}{b_3}\Big)\, .\end{eqnarray*}
Because $(1+u)^v >1+u^v$ for $u>0, v\in(0,1)\,,$ this last inequality
readily yields

\be\label{xt-mz-13} |v_1-v_2|>
\Big(\frac{b_4}{b_3}\Big)^{\frac{1}{3c'M}}\min\{d_{D_1}(v_1),d_{D_1}(v_2)\}.\ee

Without loss of
generality, we may assume that
$$\min\{d_{D_1}(v_1),d_{D_1}(v_2)\}=d_{D_1}(v_1).$$ Then we get

 \begin{lem}\label{cl3.1}  For all $z'\in\gamma'[v'_1,v'_2]$, we have $d_{D'}(z')\leq b_3^3d_{D'_1}(z')$.
\end{lem}

\bpf Suppose on the contrary that there exists some point $v_3'\in\gamma'[v'_1,v'_2]$ such
that

\be\label{eq(3-7-1-1)}d_{D'}(v'_3)>b_3^3\,d_{D'_1}(v'_3).\ee

We divide our discussions into two cases:
$$\min\{\ell(\gamma'[x',v'_3]),\ell(\gamma'[y',v'_3])\}=\ell(\gamma'[x',v'_3])$$
and
$$\min\{\ell(\gamma'[x',v'_3]),\ell(\gamma'[y',v'_3])\}=\ell(\gamma'[y',v'_3]).$$
We only need to consider the former case since the discussion for the latter one
is similar.

Obviously, Corollary \ref{cor-2} and (\ref{eq(3-7-1-1)}) imply that for all $w'\in \gamma'[x',v'_3]$,

\be\label{xt-li-1}\ell(\gamma'[w',v'_3])\leq 2b_3^2\,d_{D'_1}(v'_3)
<\frac{2}{b_3}d_{D'}(v'_3),\ee and hence for all $v'\in \gamma'[y',v'_2]$, \eqref{xt-mz-1} shows

\beq\label{hl--3} |v'_3-v'|&\leq&
|y'-v'|+|y'-x'|+|x'-v'_3| \leq
\Big(\frac{1}{12b_3^2}+2\Big)\ell(\gamma'[x',v'_3]) \\ \nonumber &\leq&
\frac{2}{b_3}\Big(\frac{1}{12b_3^2}+2\Big)\,d_{D'}(v'_3).\eeq
We know from (\ref{xt-li-1}) and \eqref{hl--3} that $x',  v'_1, v'_2\in \mathbb{B}(v'_3, \frac{1}{2}d_{D'}(v'_3))$.

\begin{figure}[!ht]
\centering
\includegraphics[width=8cm]{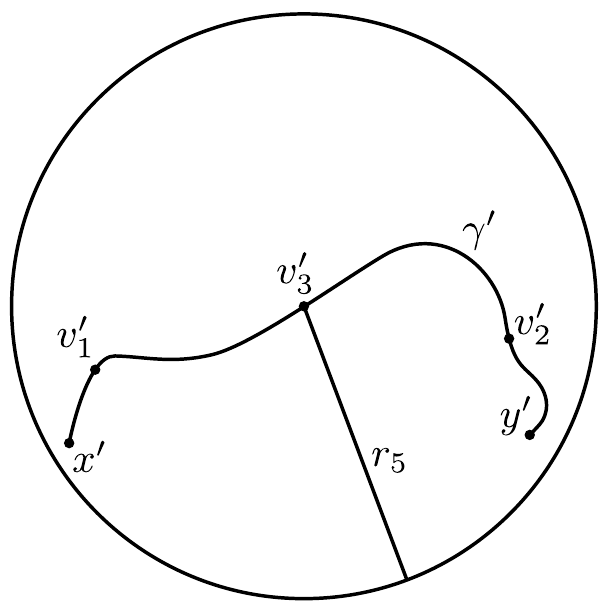}
\caption{\label{fig2} $r_5=\frac{1}{2}d_{D'}(v'_3)\,.$}
\end{figure}

Suppose $z_0\in \gamma[x, v_1]$. Then Corollary \ref{cor1} implies that

$$|v_2-v_1|\leq \rho_1d_{D_1}(v_1),$$
which contradicts \eqref{xt-mz-13}. Hence we have proved $v_1\in \gamma[x, z_0]$, and then
Theorem $H$, (\ref{ly-1}) and (\ref{xt-mz-13}) show that

\be\label{eq(3-11)}\frac{1}{\rho_1}\Big(\frac{b_4}{b_3}\Big)^{\frac{1}{3c'M}}\leq\frac{|v_1-v_2|}{|v_1-x|}\leq
\eta\Big(\frac{|v'_1-v'_2|}{|v'_1-x'|}\Big)< \eta(13b_3^2),
\ee since $|v'_1-v'_2|\leq |x'-v'_1|+|y'-v'_2|+|x'-y'|\leq (24b_3^2+1)|x'-y'|$.
This is a contradiction, which  completes the proof of Lemma \ref{cl3.1}.
\epf

\begin{lem}\label{cl3.2}  For all
$z\in\gamma[v_1,v_2]$, we have $d_{D}(z)\leq 2\tau d_{D_1}(z)$,
here and in the following, $\tau=\frac{2}{\eta^{-1}(\frac{1}{b_3^3})}$. \end{lem}

\bpf Suppose on the contrary that there is some $\zeta\in\gamma[v_1,v_2]$ such that

\be\label{xt-mz-14}d_{D}(\zeta)> 2\tau d_{D_1}(\zeta).\ee

\begin{figure}[!ht]
\centering
\includegraphics[width=12cm]{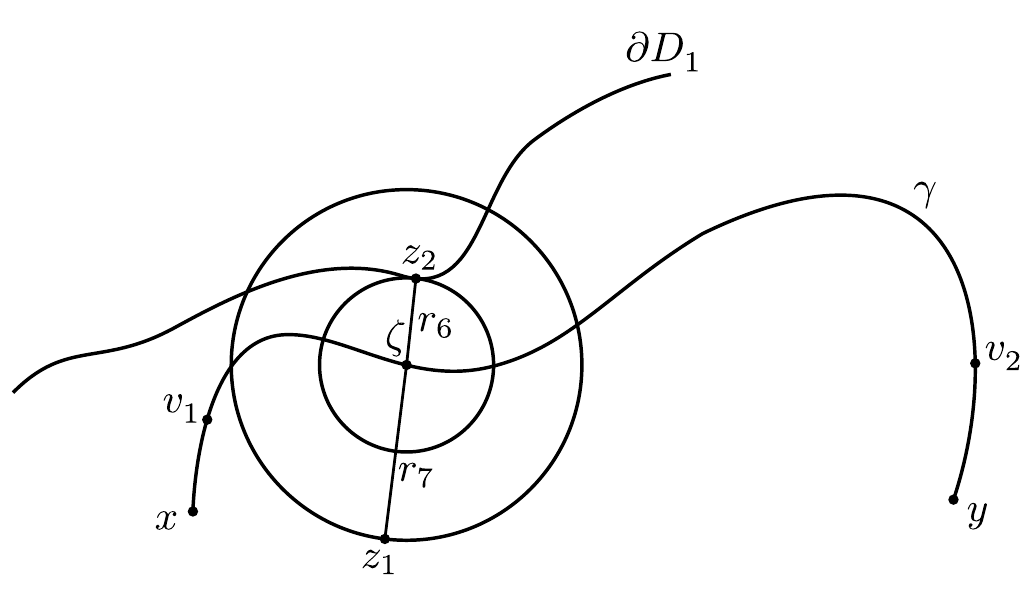}
\caption{\label{fig1} $ r_6\leq 2d_{D_1}(\zeta), r_7=\frac{1}{2}d_{D}(\zeta)\,.$}
\end{figure}

In order to get a contradiction, we let $z'_1\in
f(\mathbb{S}(\zeta,\frac{1}{2}d_{D}(\zeta)))\cap
\mathbb{B}(\zeta',d_{D'}(\zeta'))$ and $z_2\in \partial D_1$ be such that
$|z_2-\zeta|\leq 2d_{D_1}(\zeta)$. We infer from
\eqref{xt-mz-14} and Lemma \ref{cl3.1} that $z_1$, $z_2\in
\mathbb{B}(\zeta, \frac{3}{4}d_{D}(\zeta))$,
$$|z'_2-\zeta'|\geq d_{D'_1}(\zeta')\geq \frac{1}{b_3^3}d_{D'}(\zeta'),$$
$$|z'_1-\zeta'|<d_{D'}(\zeta')\;\; \mbox{and}\;\;|z_1-\zeta|=\frac{1}{2}d_{D}(\zeta)>\tau d_{D_1}(\zeta).$$
Hence
\begin{eqnarray*}\frac{1}{b_3^3}\leq \frac{|z'_2-\zeta'|}{|z'_1-\zeta'|}\leq \eta\Big(\frac{|z_2-\zeta|}{|z_1-\zeta|}\Big)<
\frac{1}{b_3^3}.\end{eqnarray*}
This obvious contradiction completes the proof of Lemma
\ref{cl3.2}.\epf

\noindent {\bf Now we are ready to finish the proof of \eqref{lem0-1}}.

\medskip

It follows from \eqref{lem-j-j} and \eqref{xt-mz-1} that
\begin{eqnarray*}
|v_1'-v_2'|\leq |x'-v_1'|+|x'-y'|+|v_2'-y'|&\leq& \Big(2+\frac{1}{12b_3^2}\Big)\min\{\ell(\gamma'(x', v_1')), \ell(\gamma'(y', v_2'))\}
\\ \nonumber &\leq & \frac{1+24b_3^2}{4}\min\{d_{D'}(v_1'), d_{D'}(v_2')
\},
\end{eqnarray*}
and then
  \eqref{Wen-1}, \eqref{Wen-2}, \eqref{lem-j-j}, \eqref{ly-0}, (\ref{xt-mz-1}), \eqref{xt-mz-12} and Lemma
\ref{cl3.2} imply that

\begin{eqnarray*}
\log\Big(1+\frac{b_4-24b_3^2}{24b_3^2}\Big) &\leq &
\log\Big(1+\frac{\ell(\gamma'[v'_1,v'_2])}{\min\{d_{D'_1}(v'_1),d_{D'_1}(v'_2)\}}\Big)
\leq \ell_{k_{D'_1}}(\gamma'[v'_1,v'_2])\\ \nonumber &\leq&
2k_{D'_1}(v'_1,v'_2)\leq 2Mk_{D_1}(v_1,v_2)+2C\\ \nonumber &\leq&
2c'M\;j_{D_1}(v_1,v_2) +2C\\ \nonumber &\leq&
2c'M\log\Big(1+\frac{2\tau
|v_1-v_2|}{\min\{d_{D}(v_1),d_{D}(v_2)\}}\Big)+2C\\ \nonumber &\leq
& 4c'\tau Mk_{D}(v_1,v_2)+2C \\ \nonumber &\leq & 4c'\tau
M^2k_{D'}(v'_1,v'_{2})+4c'\tau CM+2C\\ \nonumber &\leq&
4a'c'M^2\tau\;j_{D'}(v'_1,v'_2) +4c'CM\tau +2C \\ \nonumber &\leq&
4a'c'M^2\tau\;j_{D'_1}(v'_1,v'_2) +4c'CM\tau +2C\\ \nonumber &<&
8a'c'M^2\tau\log\frac{5+24b_3^2}{4},
\end{eqnarray*} which contradicts with $(2)$ in Corollary \ref{dear-1}.
 Then we know that \eqref{lem0-1} is true.

\subsection{The proof of Theorem \ref{thm1.1}}

The inequalities \eqref{lem-j-j}, \eqref{lem0-1} and the arbitrariness of the pair of points $x'$ and $y'$ in
$D_1'$ show that $D_1'$ is $c'$-uniform, which implies that Theorem \ref{thm1.1} is true.\qed

\bigskip
{\bf Acknowledgement.} This research was carried out when the first author was an academic visitor
 at the University of Turku and supported by the Academy of Finland grant of the second author with the
Project number 2600066611. The authors are indebted to the referee for his/her
valuable suggestions.

\end{document}